\documentclass[12pt]{amsart}

\usepackage{epsf}

\newtheorem{thm}{Theorem}

\newtheorem{lemma}[thm]{Lemma}
\newtheorem{cor}[thm]{Corollary}

\newcommand{\R}{\mbox{\bf{R}}}
\newcommand{\Z}{\mbox{\bf{Z}}}

\newcommand{\bdry}{\partial}

\def\dfn#1{{\em #1}}

\begin{document}
%%%%%%%%%%%%%%%%%%%%%%%%%%%%%%%%%%%%%%%%%%%%%%%%%%%%%%%%%
\title{On the non-existence of tight contact structures}

\author{John B. Etnyre}
\address{Stanford University, Stanford, CA 94305}
\thanks{JE supported in part by NSF Grant \# DMS-9705949}
\email{etnyre@math.stanford.edu}
\urladdr{http://math.stanford.edu/\char126 etnyre}

\author{Ko Honda}
\address{University of Georgia, Athens, GA 30602}
\email{honda@math.uga.edu}
\urladdr{http://www.math.uga.edu/\char126 honda}

\date{October 5, 1999}


\begin{abstract}
	We exhibit a 3-manifold which admits no tight contact structure.
\end{abstract}
\maketitle
%%%%%%%%%%%%%%%%%%%%%%%%%%%%%%%%%%%%%%%%%%%%%%%%%%%%%%%%%%

%----------------------------------------------------------------------------
\section{Introduction}
%----------------------------------------------------------------------------

In 1971, Martinet \cite{Martinet} showed that every 3-manifold admits a con\-tact structure. But in
the subsequent twenty years, through the work of Bennequin \cite{Bennequin} and Eliashberg 
\cite{Eliashberg92}, it became apparent that not all contact structures
are created equal in dimension 3. Specifically,
contact structures fall into one of two classes: tight or overtwisted. In this new light, what
Martinet actually showed was that every 3-manifold admits an overtwisted contact structure.
In \cite{Eliashberg89} Eliashberg classified overtwisted contact structures on closed 3-manifolds
by proving the weak homotopy equivalence of the space of overtwisted contact structures up to isotopy and
the space of 2-plane fields up to homotopy --- hence overtwisted contact structures could now be understood
via homotopy theory.  On the other
hand, it has become apparent that tight contact structures have surprising and deep 
connections to the topology of 3-manifolds, not limited to homotopy theory.
For example, Rudolph \cite{R} and Lisca and Mati\'c \cite{LM}
found connections with slice knots and slice genus, Kronheimer and Mrowka \cite{KM}
found connections with Seiberg-Witten theory, and Eliashberg and Thurston \cite{et} found connections
with foliation theory. Thus, whether or not every 3-manifold admits a tight contact structure
became a central question in 3-dimensional contact topology.

The first candidate for a 3-manifold without a tight contact structure
was the Poincar\'e homology sphere $M=\Sigma(2,3,5)$ with reverse orientation.
The difficulty of constructing a holomorphically fillable contact structure on $M$ was highlighted
in Gompf's paper \cite{Gompf}.  Subsequently Lisca \cite{L}, using techniques from
Seiberg-Witten theory, proved that $M$ has no symplectically semi-fillable contact structure.

In this paper we prove the following nonexistence result.

\begin{thm}\label{thm:main} 
	There exist no positive tight contact structures on the Poincar\'e homology sphere
	$\Sigma(2,3,5)$ with reverse orientation.
\end{thm}

This is the first example of a closed 3-manifold which does not carry a positive tight contact
structure.

\begin{cor}   Let $M$ be the Poincar\'e homology sphere with reverse orientation.
Then the connect sum $M\#\overline{M}$, where $\overline{M}$ is $M$ with the opposite
orientation, does not carry any tight contact structure, positive or negative.
\end{cor}

This follows from Theorem \ref{thm:main}, since a tight structures on a reducible manifold may be
decomposed into tight structures on its summands \cite{Colin97, Makar-Limanov??}

%----------------------------------------------------------------------------
\section{Contact topology preliminaries}\label{section:convexity}
%----------------------------------------------------------------------------

We assume the reader is familiar with the basics ideas of contact topology in
dimension 3
(see for example \cite{a, Eliashberg92}). A thorough understanding of \cite{H}
would be helpful but we include a brief summary of the ideas and terminology.
The reader might also find \cite{Giroux91, K} useful for various parts of 
this section.

In this paper the ambient manifold $M$ will be an oriented 3-manifold, and the
contact structure $\xi$ will be {\it positive}, {\it i.e.}, given by a 1-form $\alpha$ with
$\alpha\wedge d\alpha>0$.    Throughout this section we only consider
$(M,\xi)$ {\it tight}.     Also, when we refer to {\it Legendrian curves} we mean closed curves,
in contrast to {\it Legendrian arcs}.

\subsection{Convexity}
Recall an oriented embedded surface $\Sigma$ in $(M,\xi)$ is called \dfn{convex} if
there is a vector field $v$ transverse to $\Sigma$ whose flow preserves $\xi.$ Perhaps the
most important
feature of convex surfaces is the {\it dividing set}. If $\mathcal{F}$ is a singular foliation
on $\Sigma$ then a disjoint union of (properly) embedded curves $\Gamma$ is said to 
\dfn{divide} $\mathcal{F}$ if $\Gamma$ divides $\Sigma$ into two parts $\Sigma^\pm,$ $\Gamma$
is transverse to $\mathcal{F},$ and there is a vector field $X$ directing $\mathcal{F}$
and a volume form $\omega$ on $\Sigma$ such that $\pm L_X\omega>0$ on $\Sigma^\pm$ and
$X$ points transversely out of $\Sigma^+.$ If $\alpha$ is a contact 1-form for $\xi$ then
the zeros of $\alpha(v)$ provide dividing curves for the characteristic foliation $\Sigma_\xi.$
It is sometimes useful to keep in mind that the dividing curves are where $v$ is tangent to $\xi.$
An isotopy $F:\Sigma\times[0,1]\to M$ of $\Sigma$ is called \dfn{admissible} if $F(\Sigma\times\{t\})$
is transversal to $v$ for all $t.$ An important result concerning convex surfaces we will (implicitly) be
using through out this paper is:
\begin{thm}[Giroux \cite{Giroux91}]\label{giroux:main}
	Let $\Gamma$ be the dividing set for $\Sigma_\xi$ and $\mathcal{F}$ another
	singular foliation on $\Sigma$ divided by $\Gamma.$ Then there is an admissible 
	isotopy $F$ of $\Sigma$ such that $F(\Sigma\times\{0\})=\Sigma,$ 
	$F(\Sigma\times\{1\})_\xi=\mathcal{F}$ and the isotopy is fixed on $\Gamma.$
\end{thm}

Roughly speaking, this says that the dividing set $\Gamma$ dictates the geometry of $\Sigma$,
not the precise characteristic foliation.

We will let  $\#\Gamma$ denote the number of connected 
components of $\Gamma$.  If there is any ambiguity, we will also write $\Gamma_\Sigma$ instead 
of $\Gamma$ to denote the dividing set of $\Sigma$. 

\subsection{Edge-rounding} 
Let $\Sigma_1$ and $\Sigma_2$ 
be compact convex surfaces with Legendrian boundary, which 
intersect transversely along a common boundary Legendrian curve $L$.
The neighborhood of 
the common boundary Legendrian is locally isomorphic to the neighborhood $\{x^2
+ y^2  \le \varepsilon\}$ of $M = {\bf R}^2 \times ({\bf R}/{\bf Z})$ with
coordinates $(x, y, z)$ and con\-tact 1-form $\alpha = \sin(2\pi n z) dx +
\cos(2\pi n z) dy$, for some $n \in {\bf Z}^+$.  Let $A_i\subset \Sigma_i$, $i=1,2$,
be an annular collar of the boundary component $L$.
We may choose our local model so that $A_1 = \{x = 0, 0 \le
y \le \varepsilon\}$
and $A_2 = \{y = 0, 0 \le x \le \varepsilon\}$ (or the same with $A_1$ and $A_2$  switched).
Assuming the former, if we join $\Sigma_1$
and $\Sigma_2$ along 
$x = y = 0$ and round the common edge, the resulting surface is convex, and the
dividing curve 
$z = \frac{k}{2n}$ on $\Sigma_1$ will connect to the dividing curve $z = \frac{k}
{2n} - \frac{1}
{4n}$ on $\Sigma_2$, where $k = 0, \cdots, 2n - 1$.

\subsection{Bypasses}
\begin{figure}[ht]
	{\epsfysize=1.3in\centerline{\epsfbox{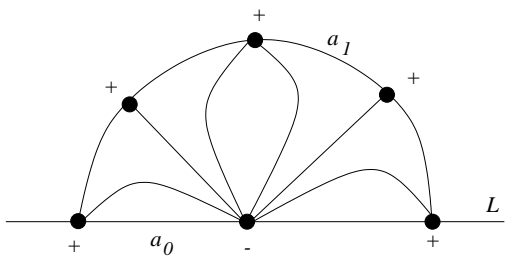}}}
	\caption{A bypass}
	\label{bypass}
\end{figure}

\begin{figure}[ht]
	{\epsfysize=1.5in\centerline{\epsfbox{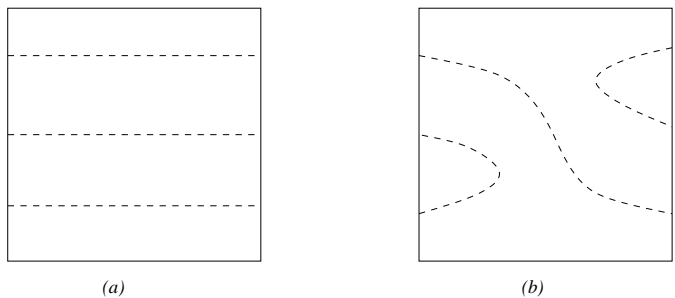}}}
	\caption{(a) Dividing curves (dashed lines) on $\Sigma$, (b) Dividing curves on $\Sigma'$}
	\label{bypassmove}
\end{figure}

We now introduce the main idea on which this work is based. Let $L$ be a Legendrian arc. A half-disk
$D$ is called a \dfn{bypass for} $L$ if $\partial D$ is Legendrian and consists of two
arcs $a_0=\partial D\cap L$ and $a_1$ such that, if we orient  $D$,
$a_0\cap a_1$ are both positive (negative)   
elliptic singular points in $D_\xi,$  there is a negative (positive) elliptic point  along $a_0$,
and the singular points along $a_1$ are all positive (negative) and alternate between elliptic and hyperbolic.
See Figure~\ref{bypass}.
We also allow for the degenerate case when $L$ is a closed curve and the    
endpoints of $a_0$ are the same ({\em i.e.}, $L=a_0$). We refer to such a bypass as a \dfn{bypass of 
degenerate type}.  Once an orientation on  
$D$ is fixed, the \dfn{sign} of the bypass 
is the sign of the elliptic point on the interior of $a_0.$
The reason for the name `bypass' is that, instead of traveling along the Legendrian arc $L$, we may
go around and drive through $L'=(L\backslash a_0)\cup a_1$   --- this has the
effect of increasing the twisting of $\xi$ along the Legendrian curve.  (Here we are using the 
convention that left twists are negative.)
One can then show:
\begin{thm}[Honda \cite{H}]\label{thm:bypass}
	Assume $M$ has a convex boundary $\Sigma$ and there exists a bypass $D$ along the
	Legendrian curve
	$L\subset \Sigma.$ Then we can find a neighborhood $N$ of $\Sigma\cup D$ with $\bdry N
	=\Sigma-\Sigma'$, and $\Gamma_{\Sigma'}$ is  related to
	$\Gamma_\Sigma$ as shown in Figure~\ref{bypassmove}.
\end{thm}

Consider a convex torus $T$. If $\xi$ is tight, then one can show that no dividing curve for
$T$  bounds a disk.  Thus $\Gamma_T$ consists of $2n$ parallel dividing curves
with $n\in\Z^+$. Assume for
now that the dividing curves are all horizontal. Using Theorem~\ref{giroux:main} we can 
isotop $T$ so that there is one closed curve of singularities in $T_\xi$ in each region of
$T\setminus \Gamma_T$ --- these are called \dfn{Legendrian divides} and are parallel to the
dividing curves. We can further assume all the other leaves of  $T_\xi$
form a 1-parameter family of parallel closed curves transverse to the Legendrian divides.
These \dfn{Legendrian ruling curves} can have any slope not equal to the slope of the dividing curves.
A convex torus $T$ is said to be in {\it standard form} if $T_\xi$  has this nongeneric form, consisting
of Legendrian divides and Legendrian ruling curves.
Suppose $n>1$ and we find a bypass for one of the ruling curves
for $T$ then we may isotop $T$ across the bypass as in Theorem~\ref{thm:bypass} and reduce the
number $n$ by one. If $n=1$ we have a configuration change as follows:
\begin{thm}[Honda \cite{H}]
	Let  $T$ be a convex torus with $\#\Gamma_T=2$ and in some basis for $T$ the slope of
	the dividing curves is 0.   If we find a bypass for a ruling curve of slope between
	$-\frac{1}{m}$ and $-\frac{1}{m+1}$, $m\in\Z,$ then after pushing $T$ across the
	bypass the new torus has two dividing curves of slope $-\frac{1}{m+1}.$  
\end{thm}

\subsection{Legendrian curves and the twisting number}
Let $\gamma$ be a Legendrian curve in $M.$ We can always find a neighborhood $N$ of $\gamma$
whose boundary is a convex torus with two dividing curves. The linking of a dividing curve
on $\partial N$ with $\gamma$ is the Thurston-Bennequin invariant of $\gamma.$  We call the
slope of these dividing curves the {\it boundary slope} of $N.$
If $\gamma$ is not null-homologous, then the framing induced on $\gamma$ relative to
some pre-assigned trivialization of the normal bundle of $\gamma$ will be called the
\dfn{twisting number} of $\gamma.$
In contrast to transverse
curves, neighborhoods of Legendrian curves have a quantifiable thickness  --- no matter
how small a (nice) neighborhood of a Legendrian curve one takes, the boundary slope of the
neighborhood is fixed.  If $N$ is the neighborhood of a Legendrian curve with twisting
number $m$ (relative to a framing that is 
already fixed), then the slopes of the dividing curves on $\bdry N$ are
$\frac{1}{m}$.  On the other hand,
any tight contact structure on the solid torus $S^1\times D^2$
with boundary slope $\frac{1}{m}$  is contact isotopic to a neighborhood of
a Legendrian curve with twisting number $m,$ see \cite{K}. 

It is easy to decrease the twisting number inside $N$ by adding `zigzags' (see below for an 
explanation of terminology).
Hence, increasing the twisting number is one way of thickening $N$. We now show how to use
bypasses to increase the twisting number of $\gamma.$
\begin{lemma}[Twist Number Lemma: Honda \cite{H}]
	Let $\gamma$ be a Legendrian curve in $M$ with a fixed framing. Let $N$ be a standard 
	neighborhood of $\gamma$ and $n$ the twisting number of $\gamma.$ 
	If there exists a bypass attached to a Legendrian ruling
	curve of $\partial N$ of slope $r$ and $\frac{1}{r}\geq n+1,$ then there exists a Legendrian
	curve with larger twisting  number isotopic to $\gamma.$
\end{lemma}
\noindent
This lemma is a very useful formulation of the observation that if one has
a bypass for a Legendrian knot then one can increase its twisting number.

\subsection{Finding bypasses}
Now that we see bypasses are useful let us consider how to find them. Let $\Sigma$ be a convex
surface with Legendrian boundary $L.$ If the Thurston-Bennequin invariant of $L$ (or
let us say twisting  number with respect to $\Sigma$) is negative
then we can arrange that all the singularities of $\Sigma_\xi$ along $L$ are (half)-elliptic. Moreover
we can assume the singular foliation $\Sigma_\xi$ is Morse-Smale.
If $n$ is the twisting number of $L$, then $\Gamma_\Sigma$ intersects $L$ at $2n$ points.
In this situation we can use
the flow of $X$ (the vector field directing $\Sigma_\xi$ in the definition of the dividing set)
to flow any dividing curve to a Legendrian arc on $\Sigma$, all of whose singularities are
of the same sign.     Now, starting with a dividing curve $c$ with endpoints along $L$
isotopic to the boundary of $\Sigma$, we
push
$c$ `back into $\Sigma$ along $X$' and eventually arrive  at a Legendrian arc
which cuts a bypass for $L$ off of $\Sigma.$ Thus to find bypasses we need only look
for these $\bdry$-compressible dividing curves.

\subsection{Layering neighborhoods of Legendrian knots}    \label{layering}
We now explore the neighborhood of a Legendrian knot in more detail. To this end
let $N$ be a standard neighborhood of a Legendrian knot $\gamma$ with a fixed framing
so that $N=S^1\times D^2$ and $\gamma$ has twisting number $n\leq -1$ in this framing. Recall
this means that $\partial N$ is convex with two dividing curves (and two Legendrian
divides) of slope $\frac{1}{n}.$ Inside of $N$ we can isotop $\gamma$ (but not
Legendrian isotop) to $\gamma'$ so as to decrease the twisting number by one.
There are actually two different ways to do this. If $\gamma$ were a knot in $\R^3$
with the standard tight contact structure then this process corresponds to stabilizing 
the knot (reducing the Thurston-Bennequin invariant) by adding zigzags,
and one can do this while increasing
or decreasing the rotation number by one.

To see how to detect this difference in
$N$ we do the following: first, fix an orientation on $N$ and orient $\gamma$ so that
$\gamma$ and $\{pt\}\times D^2$ intersect positively.
Let $N'$ be a standard neighborhood of $\gamma'$.  Then consider the
layer  $U=T^2\times [0,1]=N\backslash N'$, where we set
$T^2\times\{0\}=\bdry N'$ and $T^2\times \{1\}=\bdry N.$
We assume the Legendrian ruling curves
are vertical so that $A=S^1\times [0,1]$ is an annulus with Legendrian boundary. Note that
$S^1\times \{0\}$ intersects the dividing curves on $\partial N',$ $2n-2$ times and
$S^1\times \{1\}$ intersects the dividing curves on $\partial N,$ $2n$ times. So there
will be a bypass (in fact just one bypass) on $A$ for $S^1\times \{0\}.$ If $A$ is oriented
so that the orientation on $S^1\times \{0\}$ agrees with the one chosen above on $\gamma$,
then the sign of the bypass is what distinguishes the two possible $\gamma'$'s.
If $c$ is a curve on $S^1\times \{0\}$ whose slope is not $\frac{1}{n}$ or $\frac{1}{n+1}$
then we can make the Legendrian ruling curves on $\partial U$ parallel to $c.$ Later it
will be useful to know what the dividing curves on $c\times [0,1]$ will look like. The relative 
Euler class in the next paragraph can be useful for this purpose.

\subsection{Euler class}
Let $\xi$ be any tight structure on a toric annulus $U= T^2\times[0,1]$ (not necessarily the
same one as in the above paragraphs). 
Assume the boundary is convex
and in standard form.
Let $v_b$ be a section of $\xi\vert_{\partial U}$
which is transverse to and twists (with $\xi$) along
the Legendrian ruling curves. We also take $v_b$ to be tangent to the Legendrian divides. 
We may now form the relative Euler class $e(\xi,v_b)$ in 
$H^2(U,\partial U;\Z).$
First note that $e(\xi,v_b)$ is unchanged if we perform a $C^0$-small isotopy of $\partial U$ so as
to alter the slopes of the ruling curve. Now given an oriented curve $c$ on $T^2\times\{0\}$ we can assume
the annulus $A=c\times[0,1]$ has Legendrian boundary and is also convex. We orient $A$ so that
it induces the correct orientation on $c.$ Now we have
\begin{equation}\label{eqn:formula}
	e(\xi,v_b)(A)=\chi(A_+)-\chi(A_-),
\end{equation}
where $A_\pm$ are the positive and negative regions into which the dividing curves cut $A.$ This formula
follows from Proposition~6.6 in \cite{Kanda98} once one observes that $v_b$ may be homotoped in
$\xi\vert_{\partial U}$ so as to be tangent to and define the correct orientation on $\partial A.$

\subsection{Twisting}
We end this section by making precise the notion of {\it twisting} along $T^2\times [0,1]$.

Let $\xi_1$
be the kernel of $\alpha_1= \sin z\, dx + \cos z\, dy$ on 
$T^3=T^2\times S^1=(\R^2/\Z^2)\times S^1$ with coordinates $((x,y),z)$. This contact structure
is tight. The characteristic foliation on $T^2\times\{ p\}$ is by lines of
a fixed slope. The slope ``decreases'' as $p$ moves around $S^1$ in the positive direction.
Let $[a,b]$ be an interval in $[0,\pi]$ so that at one end $T_a=T^2\times\{ a\}$ has slope $\frac{1}{n}$
and at the other end $T_b=T^2\times\{ b\}$ has slope $\frac{1}{n+1}.$ We can $C^\infty$-perturb $T_a$ and
$T_b$ so that they are convex and each has two dividing curves of slope $\frac{1}{n}$ and
$\frac{1}{n+1},$ respectively.

Let $U=T^2\times[0,1]$ be the layer between $N'$ and $N$ as in  \ref{layering}.
It is possible to show that $T^2\times [a,b]$ is contactomorphic
to $U$ with the contact structure constructed above (see \cite{H}).   
The slopes of the characteristic foliations on $T^2\times \{pt\}$ vary from $\frac{1}{n}$ to
$\frac{1}{n+1}$ in $T^2\times [a,b]$ so they do likewise in $U.$ From this one can
conclude that if $N$ is a standard neighborhood of a Legendrian knot $\gamma$ with twist number
$n$ (with respect to some fixed framing), then for any slope in $[\frac{1}{n},0)$ one can find 
a torus around $\gamma$ whose characteristic foliation has this slope (note if $n>0$ then
$[\frac{1}{n},0)$ means $[-\infty,0)\cup[\frac{1}{n},\infty]$ where $-\infty$ is
identified with $\infty$). 

More generally, if $T^2\times [0,1]$ has boundary slopes $s_i$ for $T^2\times \{i\}$, then we may find 
convex 
tori parallel to $T^2\times \{i\}$ with any slope $s$  in $[s_1, s_0]$ (if $s_0<s_1$ then
this means $[s_1,\infty]\cup [-\infty, s_0]$).   This follows from  the classification of tight
contact structures on $T^2\times I$ --- in the proof we layer $T^2\times I$ into `thin'
toric annuli, each of which is isomorphic to $T^2\times [a,b]$ above (see \cite{H}).

%---------------------------------------------------------------------------
\section{Thickening the singular fibers}
%---------------------------------------------------------------------------

Consider the Seifert fibered space $M$ with 3 singular fibers over $S^2$.
$M$ is described by the Seifert invariants $(\frac{\beta_1}{\alpha_1},\frac{\beta_2}{\alpha_2},
\frac{\beta_3}{\alpha_3})$.  Let $V_i$, $i=1,2,3$, be the neighborhoods of the singular
fibers $F_i$, isomorphic to $S^1\times D^2$ and identify $M\setminus(\cup_i V_i)$ with
$S^1\times \Sigma_0,$ where $\Sigma_0$ is a sphere with three punctures.   Then $A_i:
\bdry V_i\rightarrow -\bdry (M\backslash V_i)$ is
given by
$A_i=\left( \begin{array}{cc} \alpha_i & \gamma_i \\ -\beta_i & \delta_i\end{array}
\right)\in SL(2,\Z)$.
We identify $\bdry V_i=\R^2/\Z^2$, by choosing $(1,0)^T$ as the meridional direction, and $(0,1)^T$ as the
longitudinal direction with respect to the product structure on $V_i$.  We identify $-\bdry (M\setminus V_i)=\R^2/\Z^2,$ by
letting
$(0,1)^T$
be the direction of an $S^1$-fiber, and $(1,0)^T$ be the direction given by $\bdry (M\setminus V_i)
\cap (\{pt\}\times \Sigma_0)$.

Let $M$ be the Poincar\'e homology sphere $\Sigma(2,3,5)$ with reverse orientation.
It is a Seifert fibered space over $S^2$ with Seifert invariants $(-\frac{1}{2},
\frac{1}{3}, \frac{1}{5})$.
In the case $V=V_1$, we choose $A_1=
\left(\begin{array}{cc} 2 & -1 \\ 1 & 0\end{array}\right)$.
Notice there we has some freedom to choose $\gamma_1$, $\delta_1$, since
we could have postmultiplied by $\left( \begin{array}{cc}
1 & m\\ 0 & 1\end{array} \right)$ if we changed our framing for $V_1$.
Similarly, let $A_2= \left( \begin{array}{cc}
3 & 1\\ -1 & 0\end{array} \right)$, and $A_3=
\left( \begin{array}{cc}
5 & 1\\ -1 & 0\end{array} \right)$.

Now let $\xi$ be a positive contact structure on $M$.
Assume $\xi$ is tight.  The goal of the
paper is to obtain a contradiction by finding overtwisted disks inside $M$.

In the first stage of the proof, we will try to thicken neighborhoods $V_i$
of the singular fibers $F_i$ (we may assume the singular fibers are Legendrian after isotopy).  
This is done by maximizing the twisting number
$m_i$ among Legendrian curves isotopic to $F_i$, subject to the condition that
all three Legendrian curves be simultaneously isotopic to $(F_1,F_2,F_3)$. 
Let $V_i$ be a standard tubular neighborhood of $F_i$ with minimal convex boundary
and boundary slope $\frac{1}{m_i}$ (note this is negative).

It is useful to note how the
dividing curves map under $A_i$.  $A_1: (m_1,1)^T \mapsto
(2m_1-1,m_1)^T$, $A_2: (m_2,1)^T\mapsto (3m_2+1, -m_2)^T$, and
$A_3: (m_3,1)^T\mapsto (5m_3+1, -m_3)^T$.  Therefore, the boundary slopes
are $\frac{m_1}{2m_1-1}$, $-\frac{m_2}{3m_2+1}$ and $-\frac{m_3}{5m_3+1}$,
when viewed on $-\bdry (M\setminus V_i)$.

\vskip.1in
\noindent
{\bf Warning:}  We will often call the same surface by different names, such as
$\bdry V_i$ and $-\bdry (M\setminus V_i)$.  Although the surfaces themselves are
the same, their  {\it identifications} with $\R^2/\Z^2$ are not.
Therefore, when we refer to slopes on $\bdry V_i$, we implicitly invoke the identification of
$\bdry V_i$ with $\R^2/\Z^2$ given in the first paragraph of this section.

\begin{lemma}  \label{increase}
We can increase $m_i$ so that $m_1=0$, $m_2=m_3=-1$, and thicken $V_i$
to $V_i'$ so that the slopes of $\bdry (M\setminus V'_i)$ are all infinite.
\end{lemma}

\proof We may modify the Legendrian rulings on both $\bdry(M\backslash V_2)$
and $\bdry (M\backslash V_3)$  so that
they are vertical.
Take a vertical  annulus $S^1\times I$ spanning from
a vertical Legendrian ruling curve on $\bdry (M\backslash V_2)$ to a vertical Legendrian
ruling curve on $\bdry (M\backslash V_3)$.
Here `vertical' means `in the direction of
the $S^1$-fibers'.   Note  $S^1\times I$ intersects the dividing curves
on $\partial V_2$ and $\partial V_3,$ $3m_2+1$ and $5m_3+1$ times respectively.

Assume $m_2, m_3\leq -1$.
If $3m_2+1\not= 5m_3 +1$,
then there exists a 
bypass  (due to the imbalance), attached along a vertical Legendrian curve of
$\bdry (M\backslash V_2)$ or $\bdry (M\backslash V_3)$.  We
transform $(0,1)^T$ via $A_i^{-1}$ to use the Twist Number Lemma.
$A_2^{-1}:(0,1)^T
\mapsto (-1,3)^T$, and $A_3^{-1}: (0,1)^T\mapsto (-1,5)^T$.  
The Legendrian rulings will have slope $-3$ and $-5$, and, therefore, 
we may increase the twisting number by 1 if the twisting number is $<-1$.  

Next, assume $3m_2+1=5m_3+1$ and there are no bypasses on the vertical
annulus.
Then we may take $S^1\times I$ to be standard, with vertical  rulings,
and parallel Legendrian divides spanning from $\bdry V_2$ to $\bdry V_3$.  Cutting along
$S^1\times I$ and rounding the corners, we obtain the torus boundary
of $M\setminus(V_2\cup V_3\cup(S^1\times I))$;  if we identify this torus with $\R^2/\Z^2$
in the same way as $\bdry (M\backslash V_1)$, then
the  boundary slope is $-\frac{m_2+m_3+1}{3m_2+1}= -{\frac{{\frac{8}{5}}m_2 +1}
{3m_2+1}}$.  When $m_2=-5$, then the slope is $-\frac{1}{2}$, and any Legendrian 
divide gives rise to an overtwisted disk --- this is because $A_1: (1,0)^T
\mapsto (2,1)^T$, which corresponds to a slope of $-\frac{1}{2}$ on
$\bdry (M\backslash V_1)$.  When $m_2<-5$, then the slope is $<-\frac{1}{2}$, which implies
that on the $S^1\times D^2=M\setminus(V_2\cup V_3\cup(S^1\times I))$, the twisting
in the radial direction is almost $\pi$.
In particular, there exists a convex torus with slope $\infty$ inside $S^1\times D^2$, and hence a
vertical Legendrian curve with zero twisting, obtained as a Legendrian divide of the
convex torus.    Any vertical annulus
taken from this vertical  Legendrian to $V_2$ (or $V_3$) will give bypasses.

We can now assume $m_2=m_3=-1$.
The boundary slopes of $-\bdry(M\setminus V_2)$ and 
$-\bdry(M\setminus V_3)$ are $-\frac{1}{2}$ and $-\frac{1}{4}$, respectively.
Again look at the vertical annulus $S^1\times I$ spanning from 
$V_2$ to $V_3$, with Legendrian boundary.  There are three possibilities:
(1) There are no bypasses along $S^1\times \{0\}$ (on the $V_2$ side)
--- in this case, we can cut
along the annulus as before, and get boundary component with boundary 
slope $-\frac{1}{2}$, contradicting tightness. 
(2) There is one bypass along $S^1\times \{0\}$ --- cutting along the annulus again,
we find that the boundary slope is $-1$.  This means that the twisting of
$T^2\times I=M \setminus(V_1\cup V_2\cup V_3\cup(S^1\times I))$ is large, and that  we have a 
convex torus in standard form with Legendrian divides of infinite slope
as before.  (3) There are two bypasses. 

In any case, there exist vertical Legendrian curves
with twisting number $0$ with respect to the framing from the fibers.
Since $A_1^{-1}: (0,1)^T \mapsto (1,2)^T$, we have bypasses for   
$V_1$ as well, and we can increase to $m_1=0$ using the Twist Number Lemma. 

Next, taking a vertical annulus from $V_2$ to one such vertical Legendrian curve,
we obtain two bypasses for $V_2$ and, similarly, we obtain four bypasses for $V_3$.
By attaching these two bypasses, we obtain a thickening $V_i'$ of $V_i$, $i=1,2,3$,
so that $\bdry (M\setminus V_i)$ has two vertical dividing curves.
 \qed

%-------------------------------------------------------------
\section{The Fibration}
%-------------------------------------------------------------

In the section we use the structure of $M \setminus V_3$ as a punctured torus bundle over
$S^1$ to complete the proof of our main theorem.  The strategy is as follows: First
we use the bundle structure to increase $m_3$ to 0 by finding a bypass along the
boundary of the punctured torus fiber. We then show how to increase $m_3$ to 1 by making the boundary
of the punctured torus fiber Legendrian with twist number 0. When $m_3=1$ the corresponding
neighborhood of $F_3$ almost contains an overtwisted disk.   By increasing it a little
further it does contain one. This will complete the proof of our
main theorem.

\subsection{Bundle structure}
We now describe the fibration as a punctured torus bundle over $S^1$.   In the previous section,
$V_i$ was the neighborhood of a Legendrian curve $F_i$ with twisting number $m_i$, and
$V'_i$ its thickening.  Write  $M\backslash (\cup_iV'_i)=S^1\times \Sigma$, where
$\Sigma$ is a 3-holed sphere.
Let $\gamma$ be an embedded arc on $\Sigma$ connecting
$V'_2$ to $V'_1$, $A=S^1\times\gamma$ an annulus connecting $V'_2$ to $V'_1$, and
$V$ a neighborhood of $A$ in $M \setminus (V'_1\cup V'_2)$.
Define $M'=V'_1\cup V'_2\cup V$, which is $M\setminus V_3$ with a $T^2\times I$ layer
removed from the boundary.
Let $D_1$ and $D_2$ be meridional disks for $V'_1$ and $V'_2,$ respectively.
The slope of $D_1$ on $\partial (M \setminus V'_1)$ is $-\frac{1}{2}$ and
the slope of $D_2$ on $\partial (M\setminus V'_2)$ is $\frac{1}{3}.$ Thus we can take two copies of $D_2$
and three copies of $D_1$ and glue them together with six copies of $V\cap\Sigma$ to
obtain a punctured torus $T$ in $M'$ with boundary on $\partial M'.$
(See Figure \ref{goodfiber}.)
Parallel copies of $T$ will fiber $M'$ (and $M \setminus V_3$)
as a punctured torus bundle over $S^1.$ 
The slope of $\partial T$ on $-\partial(M \setminus V_3)$ is $-\frac{1}{6}.$
Thus the slope on $\partial V_3$ is 1, so $T$ provides
a Seifert surface for $F_3.$

\begin{figure}[ht]
	{\epsfysize=2in\centerline{\epsfbox{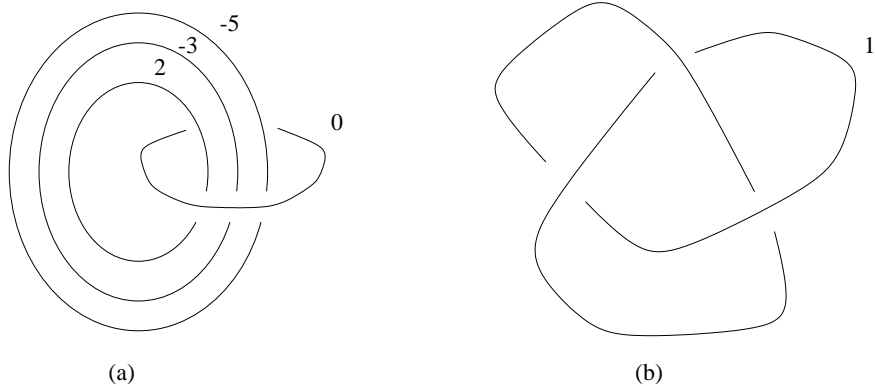}}}
	\caption{Two surgery pictures for $M$.}
	\label{fig:surgery}
\end{figure}

We can show (for example by using Kirby calculus)
that $M,$ represented in Figure~\ref{fig:surgery}(a), is diffeomorphic to
Figure~\ref{fig:surgery}(b), and
that $M \backslash V_3$ is diffeomorphic to the complement of the right-handed trefoil knot in $S^3.$
Thus we may identify the monodromy of the punctured torus bundle as 
$$\left( \begin{array}{cc}
1 & 1\\ -1 & 0\end{array} \right).$$

\subsection{Dividing curves on the punctured torus fiber}
Recall we have arranged that $m_3=-1$, so the dividing curves on $\partial V_3$
have slope $-1.$ Thus, if we make $\partial T$ a Legendrian curve on $\partial V_3$, then
$tb(\partial T)=-2.$    Here is a preliminary lemma:

\begin{lemma}
	We can always find a bypass along $\partial T$ (after possibly isotoping $T$).
\end{lemma}

\proof
We begin by showing that
either  we can find a bypass for $\partial T$ or arrange that the dividing curves
on $T$ consist of exactly two parallel arcs.
If there are no bypasses, then  the dividing set must consist of two parallel arcs
and an even number, say $2m,$ of closed parallel curves. 
Consider $M' \backslash T = T\times [0,1]$ with $T$ identified with $T\times\{0\}.$
Here we are viewing $M'$ as $M \backslash V_3$.
Let $\alpha$ be a closed Legendrian curve on $T$
parallel to the dividing curves.
Then consider $\mathcal{A}=\alpha\times[0,1]$ which we may assume is convex with Legendrian
boundary. The dividing curves do not intersect $\alpha\times\{0\}$, but they intersect $\alpha\times\{1\}$
at least $2m$ times. Thus we may find a bypass $D$ for $T\times\{1\}.$ The inner boundary $T'$ of
a neighborhood of $(T\times\{1\})\cup D$ is a convex torus with either a bypass for $\partial T'$
or two fewer dividing curves than $T.$ Repeating this argument $m$ times will result in a
bypass for $\partial T$ or a convex torus whose dividing set consist of exactly two parallel arcs.

Now suppose $\Gamma_T$  consists of two parallel arcs, since otherwise we are done.
Then $\Gamma_{T\times\{0\}}$ consists of two arcs and $\Gamma_{T\times\{1\}}$
consists of the image of these two arcs under the monodromy map for the bundle.
Let $\alpha$ be a closed curve on $T\times\{0\}$  parallel to the dividing curves, and
$\mathcal{A}= \alpha\times [0,1].$ We may arrange for
$\alpha\times \{0,1\}$ to be Legendrian and for $\mathcal{A}$ to be convex. Now the dividing
curves on $\mathcal{A}$ do not intersect $\alpha\times\{0\}$ but intersect $\alpha\times\{1\}$
at least two times. Thus we have a bypass $D$ for the dividing curves on $T\times\{1\}.$

\begin{figure}[ht]
	{\epsfysize=2in\centerline{\epsfbox{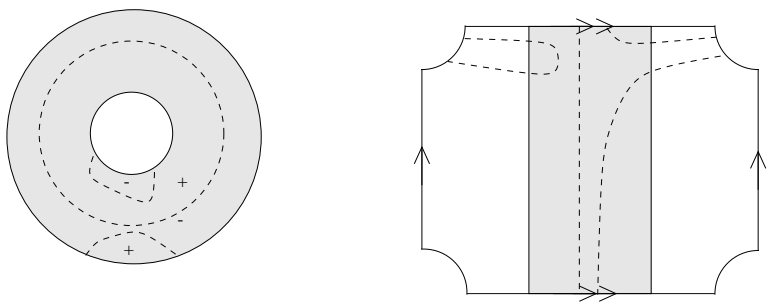}}}
	\caption{The annulus $\mathcal{A}''$ and the new $T\times\{1\}.$}
	\label{fig:annulus}
\end{figure}

Assume first that the intersection number is exactly two.
Then the annulus $\mathcal{A}$ may be split into two annuli, one of which, call it $\mathcal{A}',$
intersects $T\times\{1\}$ and contains $D$, which is of degenerate type.
 The boundary of a small neighborhood of $\mathcal{A}'$ in
$T\times(0,1)$ is an annulus $\mathcal{A}''$ and has dividing curves as shown in 
Figure~\ref{fig:annulus}.
The boundary of $\mathcal{A}''$ sits on $T\times\{1\}.$ If we cut out the annulus in $T\times\{1\}$
that cobounds a solid torus with $\mathcal{A}'',$ glue in $\mathcal{A}''$, and smooth corners,
then the dividing curves on the new $T\times\{1\}$ are shown in Figure~\ref{fig:annulus}. 
In particular 
we have a bypass for $\partial T.$

If the intersection number is greater than two, then simply attach the bypass $D$ onto
$T\times\{1\}$ --- either the new $T\times\{1\}$ has a bypass, in which case we are done,
or we can take a new $\mathcal{A}$ which has fewer (but $\geq 2$ intersections).  We
continue until the intersection number becomes two.
\qed

\subsection{Twist number increase}
We are now ready to finish the first part of the program outlined in the beginning of this 
section.
\begin{lemma}\label{lem:twist}
	We may increase
	$m_3$ to $0$ while keeping $m_1=0$ and $m_2=-1.$ Moreover, $V_2$ may be further thickened
	to $V_2''$ so that $-\partial (M\setminus V_2'')$ has dividing curves of slope $-1.$
\end{lemma}

\proof
Since we have a bypass on $\partial T$ we may apply the Twist Number Lemma to increase the twist number of
$F_3$ to $m_3=0.$ So we have two vertical dividing curves on $\partial V_3$ 
and on $\partial (M \setminus V_3)$ they have slope 0. Now repeating the argument in Lemma~\ref{increase}
we see that we can arrange $m_1=0$ and $m_2=-1.$ Note the dividing curves
on $-\partial (M\backslash V_2)$ and $\partial (M\backslash V_1)$
have slope $-\frac{1}{2}$ and 0, respectively.
Thus taking a vertical annulus between 
$\partial V_1$ and $\partial V_2$, we find a vertical bypass along $V_2.$  With
this bypass we can thicken $V_2$ to $V_2''$ whose dividing curves have slope $-1.$
\qed

Note that now $V_1$, $V_2$, $V_3$ are neighborhoods of Legendrians with $m_1=0$,
$m_2=-1$, $m_3=0$, and $V_2''$ has slope $-1$ on $-\bdry(M\setminus V_2'')$.

\subsection{Thinning before thickening}
Note if we had another bypass along $\partial T$ we could increase $m_3$ to 1 and from here we
could then find an overtwisted disk (see below). However, we know of no direct way to prove 
this bypass exists.   Therefore we use the following strategy, which can be called as `thinning
before thickening'.  Notice we have thickened $V_3$ so that $m_3=1$;  we will now backtrack
by having $V_3$ relinquish some of its thickness (we peel off a toric annulus from $V_3$ and
reattach to $V_1$ and $V_2$), and then thicken again to obtain a contradiction.

Let $D_i$ be the meridional disk to $V_i$ for $i=1,2.$ If we arrange that the  Legendrian rulings on
$\bdry V_i$ are horizontal and the $D_i$ are convex, then $\#\Gamma_{D_1}=1$ and
$\#\Gamma_{D_2}=2$. The dividing sets divide $D_1$ into two regions,
one positive and one negative, and $D_2$ into three regions --- without loss of generality we
can assume that, two regions are positive and one is negative. Note that the signs of these regions
depend on an orientation on the disks. Pick
an orientation for the fibers of the Seifert fibration, and the disks will be
oriented so that their intersection with a fiber is positive.

\begin{figure}[ht]
	{\epsfysize=2in\centerline{\epsfbox{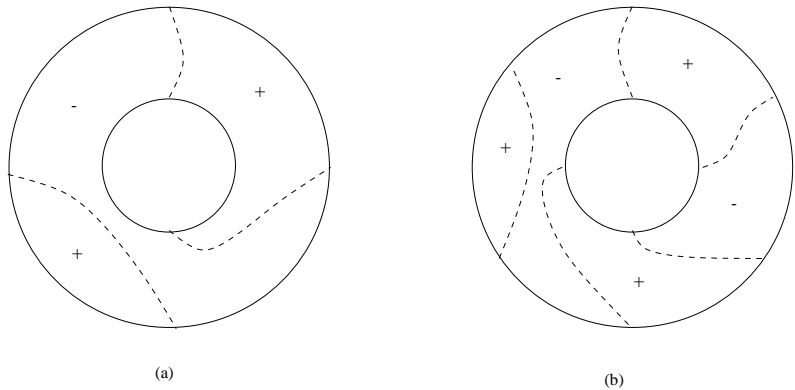}}}
	\caption{(a) Dividing curves on $c\times[0,1]$ in $V'_1\setminus V_1$, 
		(b) Dividing curves on a meridional
		disk for $V'_2$ intersected with $V'_2\setminus V_2$}
	\label{anndivides}
\end{figure}

Note we can write $V_3=C\cup U$, where $C$ is a solid torus with convex boundary and dividing
curves of slope $-1$ (on $\bdry C$), $U=T^2\times [0,1]$,
$T^2\times \{0\}=\partial C$,
and $T^2\times\{1\}$ is convex with vertical dividing curves.
We do this by stabilizing (reducing the twist number) $F_3$ with $m_3=0$ in $V_3.$
In $U$ we can find a $T^2$, say 
$T^2\times\{\frac{1}{2}\},$ with dividing curves of slope $-5$ (which correspond to vertical dividing
curves  from the point of view of $-\bdry (M\backslash V_3)$).  
By performing the correct stabilization
we can arrange that on a convex annulus in $T^2\times[\frac{1}{2},1]$ of slope $-5,$ 
with Legendrian boundary on $T^2\times \{\frac{1}{2}\}$
and $T^2\times \{1\},$  there is
a negative bypass along the boundary component on $T^2\times \{1\}.$  
We can
use this bypass to thicken $V_i$, $i=1,2$, to $V_i'$ so that $\bdry (M\setminus V_i')$ have
vertical slopes.  Using the relative Euler class
we can easily see that each meridional disk $D'_i$ in $V'_i$ has an extra positive region.
Specifically, if $e$ is the relative Euler class on $V=V'_1\setminus V_1$ (as defined in 
Section~\ref{section:convexity}) and $\mu$ and $\lambda$ are the horizontal and vertical curves
on $\partial V$ then one may easily compute that $e(\lambda\times [0,1])= 1$ and 
$e((-\mu-\lambda)\times [0,1])=0.$ Thus $e(c\times [0,1])=1$ where $c$ is the boundary of the
meridional disk in $V_1.$ From this and Equation~\ref{eqn:formula}, we can conclude that the dividing
curves are as 
in Figure~\ref{anndivides}~(a),
and hence a convex meridional disk in $V'_1$ will have one negative and two positive
regions in the complement of the dividing curves. Similarly, a convex meridional disk in
$V'_2$ will have one negative and three positive regions since a relative Euler class argument
will show its intersection with $V'_2\setminus V_2$ is shown in Figure~\ref{anndivides}~(b)

\subsection{Tight structures on $S^1\times \Sigma$}
We would now like to piece these meridional disks together to form a punctured torus
fiber for $M\setminus V'_3$,
but to do this we first need to understand the complement of the singular fibers.
Recall all the boundary slopes of $S^1\times \Sigma =M\setminus(V_1'\cup V_2'\cup V_3')$ are infinite.
Take $\Sigma=\{0\}\times\Sigma$, which we assume is convex
with Legendrian boundary.  All the boundary
components of $\Sigma$ have exactly two half-elliptic points.

\begin{lemma} Each dividing curve of $\Sigma$ must connect one boundary component
to another boundary component.
\end{lemma}

\begin{figure}[ht]
	{\epsfysize=2.6in\centerline{\epsfbox{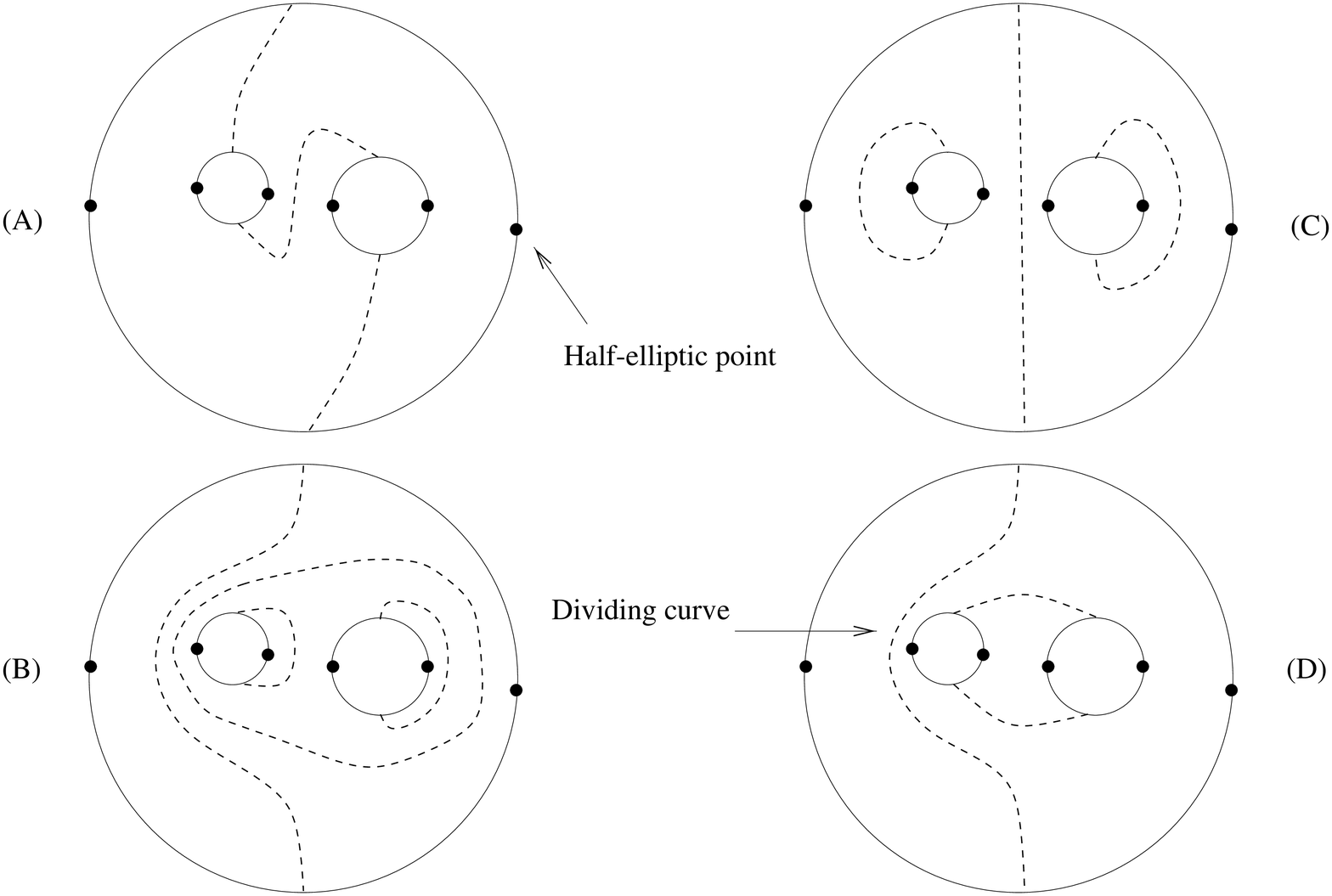}}}
	\caption{Dividing curves of $\Sigma$}
	\label{dividing}
\end{figure}

\proof  There are several possible configurations of the dividing curves.
See Figure~\ref{dividing}. 

We argue that if there is a $\bdry$-compressible
dividing arc (as in (B), (C), (D)) on $\Sigma$, then $M$ is overtwisted.   The $\bdry$-compressible
dividing arc implies the existence of a bypass (of degenerate type) along some $\bdry (M\backslash V_i')$.
Hence, there exists a layer $L_1=T^2\times I\subset S^1\times \Sigma$ with one boundary
component $\bdry (M\backslash V_i')$ and boundary slopes
$0$ and $\infty$.
However, since there are vertical Legendrian
curves with twisting number 0 outside of $V_i'\cup L_1$, we obtain another
layer $L_2=T^2\times I$, this time with boundary slopes $\infty$ and $0$.  
Therefore, $V_i'\cup L_1\cup L_2$ has too much radial twisting, and is overtwisted.
The only possible configuration without a $\bdry$-compressible dividing arc is (A).  \qed

If we take signs into consideration, there are two possible tight
structures on $S^1\times \Sigma=M\setminus(V'_1\cup V'_2\cup V'_3)$, depending on whether in
Figure~\ref{dividing}(A), the
dividing curve  from the top hits $V_2'$ or $V'_3$.

\begin{lemma}
For each of the two configurations of dividing curves on $\Sigma$,   the
tight contact structure on $S^1\times \Sigma=M\setminus(V'_1\cup V'_2\cup V'_3)$
is unique.
\end{lemma}

\proof  We cut $S^1\times \Sigma$ along $\Sigma$, round the edges, and examine
the dividing curves on the boundary of the resulting genus 2 handlebody.  We then cut along
the meridional disks of the handlebody, which we may assume have Legendrian boundary,
and eventually obtain a 3-ball.
Since each meridional
disk of the handlebody intersects the dividing set only twice, the configuration of dividing
curves on the meridional disks is unique.
Therefore the initial configuration of
dividing curves on $\Sigma$ uniquely determines the tight structure on $S^1\times \Sigma$.
\qed

\subsection{The final stretch}
\begin{figure}[ht]
	{\epsfysize=2in\centerline{\epsfbox{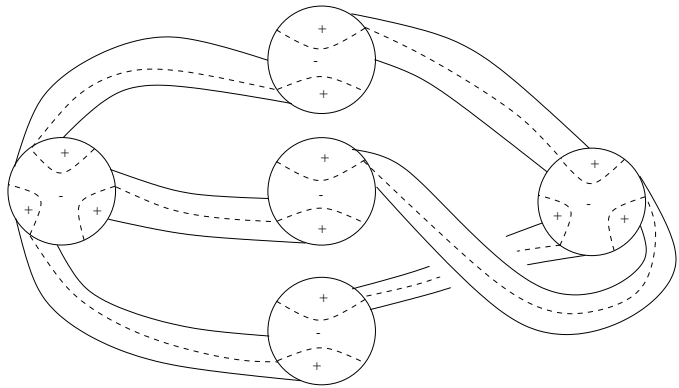}}}
	\caption{The fiber $T$ with six bypasses (note the boundaries of the disks are actually
		touching, the strips seen here just indicate how the boundaries of the disks connect)}
	\label{goodfiber}
\end{figure}

The lemma implies that the tight structure on
$S^1\times \Sigma$ is a translation-invariant tight structure on $T^2\times I$ with infinite
boundary slopes, with the standard neighborhood of a vertical Legendrian
curve with zero twisting removed. We view this $T^2\times I$  (minus $S^1\times D^2$) as the region between
$\partial V_1'$ and $\partial V_2'$ (minus $V'_3$).   We may think of the $I$-factor as being quite small, and
then isotop one of the Legendrian divides on $\partial V_1'$ to one of the Legendrian divides on
$\partial V_2'$ and finally identify small neighborhoods (in the tori) of these divides.
We may then isotop $\partial V_i'$, $i=1,2$, away from these neighborhoods so that the meridional disk $D_i'$
in $V_i'$ has Legendrian boundary. Forming the fiber $T$ in the fibration of $M\setminus V_3'$
from three copies of $D_1$ and two copies of $D_2$, we can arrange that $\partial T$ is Legendrian
and the dividing curves on $T$ are as in Figure~\ref{goodfiber}. 
Thus there are six bypasses for $\partial T.$ We can use these
to increase the twisting number of $\partial T$ to 0 ({\em i.e.}, $\partial T$ will be a
Legendrian divide on some thickening of  $\partial V_3$) which corresponds to   an
increase of $m_3$ to 1. Thus the slope of the dividing curves on $\partial (M\setminus V_3)$ is
$-\frac{1}{6}$.  Now repeat the argument in Lemma~\ref{increase} --- take a vertical
annulus from $\bdry (M\setminus V_2)$ to $\bdry (M\setminus V_3)$, and start with $m_2$
small.  Since the denominator of $-\frac{1}{6}$ is never equal to $3m_2+1$,  we eventually arrive at
$m_2=-1.$   This implies the existence of a vertical bypass for $\bdry (M\setminus V_3)$,
with which we can increase the slope of the dividing curves on $\partial (M\setminus V_3)$ to
$-\frac{1}{5}$. Thus the Legendrian divide on our thickened solid torus bounds a meridional
disk in $V_3$, so we have found an overtwisted disk.  This completes the proof of the main theorem.     \qed

\end{document}